\documentclass{gtmon_a}
\pdfoutput=1

\usepackage{stmaryrd}

\proceedingstitle{Proceedings of the School and Conference in Algebraic
Topology (The Vietnam National University, Hanoi, 9-20 August 2004)}
\conferencestart{9 August 2004}
\conferenceend{20 August 2004}
\conferencename{School and Conference in Algebraic Topology}
\conferencelocation{Vietnam National University, Hanoi, Vietnam}

\editor{John Hubbuck}
\givenname{John}
\surname{Hubbuck}

\editor{Nguy\~\ecircumflex{}n H V H\uhorn{}ng}
\givenname{H\uhorn{}ng}
\surname{Nguy\~\ecircumflex{}n}

\editor{Lionel Schwartz}
\givenname{Lionel}
\surname{Schwartz}

\title{On the construction of permutation complexes for profinite groups}

\author{Peter Symonds}
\givenname{Peter}
\surname{Symonds}
\address{School of Mathematics\\
     University of Manchester\\\newline
     PO Box 88 \\
     Manchester M60 1QD\\
     United Kingdom}
\email{Peter.Symonds@manchester.ac.uk}
\urladdr{}

\volumenumber{11}
\issuenumber{}
\publicationyear{2007}
\papernumber{17}
\startpage{369}
\endpage{378}

\doi{}
\MR{}
\Zbl{}

\arxivreference{}



\received{12 January 2005}
\revised{20 May 2005}
\accepted{1 June 2005}
\published{14 November 2007}
\publishedonline{14 November 2007}
\proposed{}
\seconded{}
\corresponding{}
\editor{}
\version{}

\makeatletter
\def\cnewtheorem#1[#2]#3{\newtheorem{#1}{#3}[section]
\expandafter\let\csname c@#1\endcsname\c@theorem}

\theoremstyle{plain}

\cnewtheorem{lemma}[theorem]{Lemma}

\cnewtheorem{proposition}[theorem]{Proposition}

\cnewtheorem{corollary}[theorem]{Corollary}

\theoremstyle{definition}

\cnewtheorem{definition}[theorem]{Definition}

\makeatother  

\theoremstyle{remark}
\newtheorem*{remark}{Remark}

 \DeclareMathOperator{\Ext}{Ext}
\DeclareMathOperator{\Hom}{Hom} \DeclareMathOperator{\res}{res}
\DeclareMathOperator{\Ker}{Ker} 
 
\DeclareMathOperator{\Res}{Res} \DeclareMathOperator{\Ind}{Ind}
 \DeclareMathOperator{\End}{End}

 \DeclareMathOperator{\Ima}{Im}
\DeclareMathOperator{\invlim}{\varprojlim}
\DeclareMathOperator{\dirlim}{\varinjlim}

 \DeclareMathOperator{\St}{St}
 
\DeclareMathOperator{\Gal}{Gal} \DeclareMathOperator{\vcd}{vcd}

\begin{document}

\begin{asciiabstract} 
Goerss, Henn, Mahowald and Rezk construct a complex of permutation 
modules for the Morava stabilizer group G_2 at the prime 
3. We describe how this can be done using techniques from homological 
algebra.
\end{asciiabstract}

\begin{htmlabstract}
Goerss, Henn, Mahowald and Rezk construct a complex of permutation
modules for the Morava stabilizer group <b>G</b><sub>2</sub> at the prime
3. We describe how this can be done using techniques from homological
algebra.
\end{htmlabstract}

\begin{abstract} 
Goerss, Henn, Mahowald and Rezk construct a complex of permutation 
modules for the Morava stabilizer group $\mathbb G_2$ at the prime 
3. We describe how this can be done using techniques from homological 
algebra.
\end{abstract}

\maketitle

\section{Introduction}
\label{intro}

In \cite{ghmr}, Goerss, Henn, Mahowald and Rezk consider the special
extended Morava stabilizer group $\mathbb G^1_2 = \mathbb S^1_2
\rtimes \Gal$ at the prime 3 and construct an exact sequence of
compact modules
$$0 \rightarrow \Ind ^{\mathbb G^1_2}_{G_{24}} \what 
{\mathbb Z} _3 \rightarrow \Ind ^{\mathbb G^1_2}_{SD_{16}} 
\what {\mathbb Z}_3(\chi) \rightarrow \Ind ^{\mathbb G^1_2}_{SD_{16}} 
\what {\mathbb Z}_3(\chi) \rightarrow  \Ind ^{\mathbb G^1_2}_{G_{24}} 
\what {\mathbb Z} _3 \rightarrow \what {\mathbb Z} _3 \rightarrow 0,
$$
where $G_{24}$ is a subgroup of order 24 etc, and $\what {\mathbb Z}
_3 (\chi )$ is a copy of $\what {\mathbb Z}_3$ on which $SD_{16}$
acts via a character $\chi \co {SD_{16}} \rightarrow \{ \pm 1 \}$.
They then use this to construct a certain tower of spectra.

The aim of this note is to show how methods from the homological
algebra and representation theory of these groups can help in the
algebraic part of this construction.

\section{Background}

Let $G$ be a profinite group and let $R$ be a complete noetherian
local ring with finite residue class field $k$ of characteristic
$p$. For example, $R$ could be the $p$--adic integers.

We work in the category of compact $R\llbracket G\rrbracket $--modules, 
$\mathcal C_R(G)$, (see Symonds
\cite{s2} for definitions, properties and more
references).

The next result is basic, but does not seem to have appeared in the
literature.

\begin{proposition}
\label{ks} If $G$ is a virtual pro--$p$--group then the Krull--Schmidt
property holds for (topologically) finitely generated modules in
$\mathcal C_R(G)$, ie every such module can be expressed as a
finite sum of indecomposable modules and this decomposition is
essentially unique in the sense that the multiplicity of each
isomorphism type is the same in any such decomposition.
\end{proposition}

\begin{proof}
Let $H \unlhd _o G$ be an open normal pro--$p$ subgroup. If $M$ is a
finitely generated $R\llbracket G\rrbracket $--module then $k \otimes _{R\llbracket H\rrbracket }M$ is
finite dimensional and we can decompose $M$ as a finite sum of
indecomposable modules using induction on 
$\dim _k k  \otimes_{R\llbracket H\rrbracket }M$.

For this to work we need to know that our induction starts, that is
that if $M \ne 0$ then $k  \otimes_{R\llbracket H\rrbracket } M \ne 0$. Let $M'$ be a
finite quotient of $M$ as an $H$--module; there is a surjection 
$k \otimes _{R\llbracket H\rrbracket } M \rightarrow k  \otimes _{R\llbracket H\rrbracket } M'$. The action
of $H$ on $M'$ factors through that of a finite $p$--group $P$, and
in this case it is well known that $k  \otimes _{R[P]} M' \ne 0$.

All we need to do now is to show that the endomorphism ring of 
a finitely generated indecomposable module is local, because then 
the uniqueness of decomposition follows formally 
(see, for example, Benson \cite[1.4.3]{benson1}).

The proof is just a variant of the one for finite groups (see
\cite[1.9]{benson1}). Let $J$ be the Jacobson radical of $R\llbracket G\rrbracket $.
For any open normal subgroup $N$ of $G$ let $I_N$ denote the
augmentation ideal of $R\llbracket N\rrbracket $. Given an endomorphism $f$ of $M \in
\mathcal C _R(G)$ we set $\Ima (f^ \infty) = \cap _{n=1}^\infty \Ima
(f^n)$ and $\Ker (f^\infty) = \{ x \in M | \forall N \unlhd _o G \,
\forall n \ge 0 \, \exists m \ge 0 \text{ such that }\  f^m (x)
\in J^nM + I_NM \}$.

For each open normal subgroup $N \unlhd _o G$ define $M_N = R
\otimes _{R\llbracket N\rrbracket } M \cong M/I_NM \in \mathcal C_R(G/N)$. Then $M
\cong \invlim M_N$. Since $M$ is finitely generated, $M_N$ is too.
Now $f$ induces an endomorphism $f_N$ of $M_N$. Define $\Ima (f_N^
\infty) = \cap _{n=1}^\infty \Ima (f_N^n)$ and $\Ker (f^\infty_N) =
\{ x \in M_N |  \forall n \ge 0 \, \exists m \ge 0 \text{ such
that } f^m (x) \in J^nM_N  \}$. From the finite group case of
Fitting's Lemma we know that $M_N = \Ima (f_N^\infty) \oplus \Ker
(f_N^\infty)$.

But $\Ima (f^\infty) \cong \invlim \Ima (f_N^\infty)$ and $\Ker
(f^\infty) \cong \invlim \Ker (f_N^\infty)$. Hence $M = \Ima
(f^\infty) \oplus \Ker (f^\infty)$.

Suppose that $M$ is indecomposable and let $I$ be a maximal left
ideal in $\End _{\mathcal C_R(G)}(M)$ and let $a$ be an endomorphism
not in $I$. Then $1 = ba + f$ for some $b \in \End _{\mathcal
C_R(G)}(M)$ and $f \in I$. But $f$ is not an isomorphism, so $M =
\Ker (f^\infty)$ and $\Ima (f^\infty )=0$.

Now $(1+f+\cdots +f^{n-1})ba= 1-f^n$. Let $N \unlhd _o G$ be some
arbitrary open normal pro-$p$ subgroup. Since $M$ is finitely
generated, for sufficiently large $n$ we have $f^n(M) \subseteq JM +
I_NM \subseteq JM$. Thus $1-f^n$ is onto, by the profinite version
of Nakayama's Lemma \cite[1.4]{brumer}. Also if $(1-f^n)(x)=0$ then
$x \in \Ima (f^\infty)=0$, so $1-f^n$ is injective. Thus $1-f^n$ is
an isomorphism and $a$ has a left inverse, $c$ say.

But $c_N$ must also be a right inverse to $a_N$ on each $M_N$, so
$c$ is also a right inverse and $a$ is an isomorphism, as required.
\end{proof}

Projective covers exist in $\mathcal C_R(G)$ (Symonds
\cite{s1}), thus so do
minimal projective resolutions.

If $S$ is a simple module, let $P_S$ denote the projective cover of
$S$. The $P_S$ are precisely the indecomposable projective modules,
and any other projective is a product of them.

If there is an open normal pro--$p$ subgroup $H \unlhd _o G$, then
any simple module for $R\llbracket G\rrbracket $ is the inflation of one for $k[G/H]$
so, in particular, there are only finitely many simple modules up to
isomorphism.

The next result is well known for finite groups.

\begin{proposition}
\label{count} Suppose that $M \in \mathcal C_R(G)$ is projective
over $R$ and let
$$\cdots \rightarrow P_r \rightarrow
  \cdots \rightarrow P_1 \rightarrow P_0 \rightarrow M$$
be the minimal
projective resolution of $M$. If $S$ is a simple module then the
multiplicity of $P_S$ in $P_r$ is equal to $\dim _{\End (S)}\Ext
^r_{R\llbracket G\rrbracket }(M,S)= \dim _{\End (S)} H^r(G, (k \otimes _R M)^* \otimes
_R S)$.
\end{proposition}

Here $S^*$ denotes the dual over $k$, or rather the contragredient.

(If $k$ is a splitting field for $G/H$, where $H<G$ is open, normal
and pro-$p$, then $\End (S) \cong k$.)
\begin{proof}(cf Symonds--Weigel \cite{sw})\qua
The multiplicity of $P_S$ in $P_r$ is
$$\dim _{\End (S)} \Hom _{R\llbracket G\rrbracket }(P_r,S).$$
The fact that the 
projective resolution is minimal implies that the differentials in 
the complex $\Hom _{R\llbracket G\rrbracket }(P_{\bullet},S)$ are zero.

Combining these facts, we find that the multiplicity is
$\dim _{\End (S)} \Ext ^r _{R\llbracket G\rrbracket }(M,S).$
But 
$\Ext^r_{R\llbracket G\rrbracket }(M,S) \cong \Ext ^r_{R\llbracket G\rrbracket }(R,\Hom _R(M,S))$  
(see eg \cite[3.1.8]{benson1}) and 
$\Hom _R(M,S) \cong (k \otimes _R M)^*\otimes _RS.$
\end{proof}

From now on we assume that $G$ is of finite virtual cohomological
dimension over $R$. The definition of Tate--Farrell cohomology
appears in Scheiderer \cite{scheiderer} for discrete coefficients and in
Symonds \cite{s2} for compact ones, as does the next result. (See
Brown \cite{brown} for its basic properties in the case of an abstract
group.)

\begin{proposition}
For $M$ in $\mathcal C_R(G)$ or $\mathcal D_R(G)$, the Tate--Farrell
cohomology $\what{H}^*(G,M)$ is isomorphic to the equivariant
Tate--Farrell cohomology of the Quillen complex of $G$ with
coefficients in $M$.
\end{proposition}

\begin{corollary}
\label{rank1} If $G$ has $p$--rank 1 (ie no subgroups isomorphic
to $\mathbb Z/p \times \mathbb Z/p$) and only finitely many
conjugacy classes of subgroups isomorphic to $\mathbb Z /p$ with
representatives $C_1, \ldots , C_n$ then $\what H^*(G,M) \cong \oplus
_{i=1}^n \what H^*(N_G(C_i),M)$ for any $M$ in $\mathcal C_R(G)$ or
$\mathcal D_R(G)$.
\end{corollary}

A similar result for $M=k$ also appears in Henn \cite{henn}.

For $M,N \in \mathcal C_R(G)$ we can also define Tate--Farrell
$\Ext$ groups $\smash{\widehat \Ext}^*_G(M,N)$. This allows us to define the
stable category $\St _R(G)$ to have the same objects as $\mathcal
C_R(G)$ but morphism groups $\smash{\widehat \Ext}^0_G(M,N)$. We write $\simeq$
for isomorphism in the stable category.

There is another description. We define the Heller translate
$\Omega$ on $\mathcal C_R(G)$ by the short exact sequence $\Omega M
\rightarrow P_M \rightarrow M$, where $P_M$ denotes the projective
cover of $M$. We also define $\underline{\Hom} _G(M,N)$ to be the
quotient of $\Hom _{\mathcal C _R(G) }(M,N)$ by the submodule of all
homomorphisms that factor through a projective module. Then $\smash{\widehat
\Ext}^r_G(M,N) \cong \dirlim _i \underline{\Hom}_G(\Omega
^{r+i}M,\Omega ^iN)$. In fact we only need to take $i \ge \vcd G$.

For the basic properties of the stable category see 
Benson \cite{benson1} for finite groups and 
\cite{benson2} for infinite abstract groups.
In particular, it is a triangulated category with the inverse of
$\Omega$ as translation and the exact triangles coming from short
exact sequences in $\mathcal C_R(G)$.

The next statement is basic to our approach, although it is just a
corollary of Yoneda's Lemma.

\begin{lemma}
\label{yon} If the homomorphism $f \co A \rightarrow B$ induces an
isomorphism
$$f^* \co \smash{\widehat\Ext}^0_G(B,M) \rightarrow
  \smash{\widehat\Ext}_G^0(A,M)$$
for all $M \in \mathcal C_R(G)$ then $f$ is an
isomorphism in the stable category.
\end{lemma}

\begin{definition}
A module $M \in \mathcal C_R(G)$ is \emph{cofibrant\/} if it is
projective on restriction to some open subgroup of $G$.
\end{definition}

In fact, if $M$ is cofibrant then it is projective on restriction to
any $p$--torsion free subgroup.

Notice that $\Omega ^i M$ is always cofibrant if $i \ge \vcd G$. If
$M$ and $N$ are cofibrant then $\smash{\widehat\Ext}^0_G(M,N) \cong
\underline{\Hom}_G(M,N)$.

The definition is taken from \cite{benson2}, as is the next lemma.
As the terminology suggests, this is part of a the structure of a
closed model category, but we do not need that here.

\begin{lemma}
\label{cof} If $M \simeq N$ in $\St _R(G)$ and $M$ and $N$ are
cofibrant then there exist projective modules $P$ and $Q$ such that
$M \oplus P \cong N \oplus Q$ in $\mathcal C_R(G)$. If $M$ and $N$
are finitely generated then $P$ and $Q$ can be chosen to be finitely
generated.
\end{lemma}

\begin{proof}
Let $H \unlhd _oG$ be open normal of finite cohomological dimension.
The inclusion of the fixed points induces a map $R \rightarrow
R[G/H]$, which is split over $H$. This induces a map $M \rightarrow
R[G/H] \otimes M \cong \Ind ^G_HM$, which is also split over $H$ and
where $Q = \Ind ^G_HM$ is projective, and finitely generated if $M$
is.

Consider the map $M \rightarrow Q \oplus N$, where the first
component is the map constructed above and the second is a stable
isomorphism. This map is split over $H$, so the cokernel, call it
$P$, is cofibrant, and finitely generated if $M$ and $N$ are.

The long exact sequence for $\smash{\widehat \Ext}^*_G(P,-)$ tells us that
$0=\smash{\widehat \Ext}^0_G(P,P) \cong \underline{\Hom}_G(P,P)$, so $P$ is
projective and the short exact sequence splits.
\end{proof}

\section{The calculation}

We set $R= \what {\mathbb Z}_3$, $k=\mathbb F_3$. The Morava
stabilizer group $\mathbb S_2$ at the prime 3 can be split as a
product $\mathbb S^1_2 \times \what {\mathbb Z}_3$, where $\mathbb
S^1_2$ is the kernel of the reduced norm. There is a natural action
of the Galois group $\Gal = \Gal (\mathbb F_9 / \mathbb F_3)$, and
we will consider the special extended Morava stabilizer group
$\mathbb G^1_2= \mathbb S^1_2 \rtimes \Gal$.

Let $S^1_2$ be the Sylow 3--subgroup of $\mathbb S^1_2$. It is normal
in $\mathbb G^1_2$ and $\mathbb G^1_2 = S^1_2 \rtimes SD_{16}$,
where $SD_{16}$ is a subgroup isomorphic to the special dihedral
group of order 16. In fact, if $\phi$ denotes the generator of
$\Gal$ (of order 2) there is an element $\omega \in S^1_2$ of order
8 such that $SD_{16}$ is generated by $\phi$ and $\omega$. There is
just one finite 3-subgroup, up to conjugation. It is cyclic of order
3 and we denote it by $C_3$. It is contained in a subgroup $G_{24}$
of order 24, but there is no subgroup of order 48. We can, however,
choose conjugacy class representatives so that $SD_{16} \cap G_{24}
= Q_8$, a quaternion group of order 8 generated by $\omega \phi$,
which commutes with $C_3$, and $\omega ^2$, which does not. We refer
to \cite{ghmr} for the details.

As a consequence, the simple modules in $\mathcal C _{\what {\mathbb
Z}_3}(G)$ correspond to the simple modules for $SD_{16}$ over
${\mathbb F_3}$. In particular there is a character $\chi$
corresponding to the map $SD_{16} \rightarrow SD_{16}/Q_8 \cong \{
\pm 1 \}$, so $\chi (\phi ) =  \chi (\omega ) =-1$. Define a module
$N_1$ by $$
0 \rightarrow N_1  \rightarrow \Ind ^{\mathbb
G^1_2}_{G_{24}} {\what {\mathbb Z}_3} \rightarrow {\what {\mathbb
Z}_3} \rightarrow 0,
$$ 
where the right hand arrow is the natural
augmentation.

Let $S$ be a simple module and apply $\Ext ^*_{\mathbb G^1_2}(-,S)$.
We obtain the long exact sequence $$\cdots \rightarrow \Ext
_{\mathbb G^1_2}^*({\what {\mathbb Z}_3},S) \rightarrow \Ext
_{\mathbb G^1_2} ^*(\Ind ^{\mathbb G^1_2}_{G_{24}}{\what {\mathbb
Z}_3},S) \rightarrow \Ext _{\mathbb G^1_2}^*(N_1,S) \rightarrow
\cdots.$$
The arrow on the left is just $\smash{H^*({\mathbb G^1_2},S)
\stackrel{\res }{\longrightarrow} H^*(G_{24},S)}$, which is equivalent
to
$$H^*(S^1_2,S)^{SD_{16}} \stackrel{\res}{\longrightarrow}
  H^*(C_3,S)^{C_8} \quad\text{or}\quad
  (H^*(S^1_2)\otimes S)^{SD_{16}}
  \stackrel{\res}{\longrightarrow} (H^*(C_3)\otimes S)^{C_8}$$
or, more
naturally,
$$(H^*(S^1_2)\otimes S)^{SD_{16}} \stackrel{\res}{\longrightarrow} ((H^*(C_3)
  \oplus H^*(C'_3)) \otimes S)^{SD_{16}},$$
where $C'_3$ is the conjugate of $C_3$ by $\omega$.
(Where no coefficients for the cohomology are indicated they are
just ${\mathbb F_3}$.)

Now, for any finite ${\mathbb F_3}SD_{16}$--module $A$, the number
$$
\dim _{\End (S)}(A \otimes S)^{SD_{16}} \cong \dim _{\End(S^*)}
\Hom _{SD_{16}}(S^*,A)
$$
is just the multiplicity of the dual $S^*$
as a summand of $A$ ($A$ is completely reducible). So we are just
decomposing the ${\mathbb F_3}SD_{16}$--modules and identifying the
map $\rho \co H^*(S^1_2) \stackrel{\res}{\rightarrow} H^*(C_3)
\oplus H^*(C'_3)$. But this factors as $$H^*(S^1_2)
\stackrel{\res}{\rightarrow} H^*(C_{S^1_2}(C_3)) \oplus
H^*(C_{S^1_2}(C'_3)) \stackrel{\res}{\rightarrow} H^*(C_3) \oplus
H^*(C'_3).$$ A standard calculation \cite{ghmr,henn,s1} shows
that $C_{S^1_2}(C_3) \cong \what {\mathbb Z}_3 \times C_3$. Its
cohomology is just
$$H^*(C_{S^1_2}(C_3)) \cong H^*(\mathbb Z_3)
\otimes H^*(C_3) \cong E(a_1) \otimes ({\mathbb F_3}[y_1] \otimes
E(x_1)) \cong {\mathbb F_3}[y_1] \otimes E(x_1,a_1),$$
where $E$ denotes an exterior algebra, $a_1,x_1$ are in degree 1 and
$y_1$ is in degree 2. The restriction to $C_3$ just kills $a_1$. For
$\smash{C_{S^1_2}(C'_3)}$ the result is similar, but we use the subscript
2 for the generators, which we take to be the images of those in the
first case under conjugation by $\omega$.

Henn \cite{henn} shows that the first of the maps above is
injective. Its image is generated as an algebra by
$x_1,x_2,y_1,y_2,(x_1a_1-x_2a_2),y_1a_1,y_2a_2$.  The action of
$SD_{16}$ can now be calculated and is given in \cite{ghmr}:
\begin{align*}
\omega _* (x_i) & =-(-1)^ix_{i+1}, & 
\omega _* (y_i) & =-(-1)^iy_{i+1}, &
\omega _* (a_i) & =-(-1)^ia_{i+1}, \\
\phi _* (x_i) & =-x_{i+1},  &
\phi _* (y_i) & =-y_{i+1}, &
\phi _* (a_i) & =-a_{i+1},
\end{align*} 
(where the subscripts are taken modulo 2).

The map $\rho$ is also explicitly calculated in 
Gorbounov--Siegel--Symonds \cite{gss}.

From this we can read off that $\rho$ is surjective, except in
degree 0, where the cokernel is ${\mathbb F_3}(\chi)$ as an
$SD_{16}$--module. It is also injective in degrees 0 and 1. In degree
2 the kernel is generated by $x_1a_1-x_2a_2$, which gives a copy of
${\mathbb F_3}(\chi)$ again. In degree 3 the kernel is generated by
$y_1a_1$ and $y_2a_2$, so consists of two simples: one trivial
generated by $y_1a_1+y_2a_2$ and a copy of ${\mathbb F_3}(\chi)$
generated by $y_1a_1-y_2a_2$.

Thus the minimal projective resolution of $N_1$ starts 
$$ 
\cdots
\rightarrow P_{\mathbb F_3} \oplus P_{{\mathbb F_3}(\chi)}
\rightarrow P_{{\mathbb F_3}(\chi)} \rightarrow P_{{\mathbb
F_3}(\chi)} \rightarrow N_1 \rightarrow 0.$$ 
Now 
$P_{{\mathbb F_3}(\chi)} \cong \Ind ^G_{SD_{16}}
{\what {\mathbb Z}_3}(\chi)$,
because the latter is projective and, for any simple $S$,
$$\Hom_G(\Ind ^G_{SD_{16}}{\what {\mathbb Z}_3}(\chi),S) \cong \Hom
_{SD_{16}}({\what {\mathbb Z}_3}(\chi),S),$$
which is non-zero only
for $S\cong {\mathbb F_3}(\chi)$ and then it has dimension $1$. So
if we define $N_3=\Omega ^2 N_1$ we have an exact sequence
$$
0 \rightarrow N_3 \rightarrow \Ind ^G_{SD_{16}}
{\what {\mathbb Z}_3}(\chi) \rightarrow \Ind ^G_{SD_{16}}
{\what {\mathbb Z}_3}(\chi) \rightarrow N_1 \rightarrow 0,$$
where $N_3$ has projective cover $P_{\mathbb F_3} \oplus P_{{\mathbb
F_3}(\chi)}$.

If we work stably we can obtain $\Omega ^2 N_1$ another way. Recall
that $C_3$ is the only cyclic subgroup of order 3 in $\mathbb G^1_2$
up to conjugacy. Write $N=\smash{N_{\mathbb G^1_2}(C_3)}$; because $Q_8$
normalizes $C_3$ it also normalizes $\smash{C_{S^1_2}(C_3)}$, and since the
centralizer can be of index at most 2 in the normalizer we see that
$N \cong C_3 \times \what {\mathbb Z}_3 \rtimes Q_8$.

From \fullref{rank1} we see that
$$\res \co \what H^*({\mathbb G^1_2},M) \rightarrow \what H^*(N,M)$$
is an isomorphism, or equivalently that
the augmentation map
$\smash{\epsilon \co \Ind ^{\mathbb G^1_2}_N{\what
{\mathbb Z}_3} \rightarrow {\what {\mathbb Z}_3}}$
induces an
isomorphism $\smash{\widehat \Ext^*_{\mathbb G^1_2}({\what {\mathbb
Z}_3},M)}
\rightarrow \smash{\widehat \Ext _{\mathbb G^1_2}^*(\Ind ^{\mathbb
G^1_2}_N{\what {\mathbb Z}_3},M)}$, for any $M \in \mathcal C_{\what
{\mathbb Z}_3}(G)$. It follows from \fullref{yon} that $\epsilon$ is a
stable isomorphism.

So stably our complex starts
$$\Ind ^{\mathbb G^1_2}_{G_{24}}{\what
{\mathbb Z}_3} \rightarrow \Ind ^{\mathbb G^1_2}_N{\what {\mathbb
Z}_3},$$ which is $\Ind ^{\mathbb G^1_2}_N$ applied to the natural
augmentation map $\Ind ^N_{G_{24}}{\what {\mathbb Z}_3} \rightarrow
{\what {\mathbb Z}_3}$ over $N$.

But the subgroup $D < G_{24}$ generated by $C_3$ and $\omega \phi$
is normal in $N$, so $N$ acts on $\Ind ^N_{G_{24}}{\what {\mathbb
Z}_3}$ via its image $N/D \cong \what {\mathbb Z}_3 \rtimes C_2$, the
infinite virtually 3-adic dihedral group, so we can resolve to
obtain 
\begin{equation}
0 \rightarrow \Ind ^N_{G_{24}}{\what {\mathbb
Z}_3}(\theta) \rightarrow \Ind ^N_{G_{24}}{\what {\mathbb Z}_3}
\rightarrow {\what {\mathbb Z}_3} \rightarrow
0,
\label{eqn:eq1}
\end{equation} where $\theta \co Q_8 \rightarrow \{ \pm
1 \}$ is the character with $\theta (\omega \phi )=1$ and $\theta
(\omega ^2)=-1$.

This can be seen systematically using cohomology, as before. More
explicitly, the non-zero map on the left is determined by $1 \otimes
\theta \mapsto (g-g^{-1}) \otimes 1$, where $g$ is a generator of
the group  $\what {\mathbb Z}_3$ and $\theta$ is considered as a
basis element of ${\what {\mathbb Z}_3}(\theta)$. The sequence is
exact because on restriction to $\what {\mathbb Z}_3$ it is just a
variation on the standard projective resolution for $\what {\mathbb
Z}_3$.

Similarly, since $G_{24}$ has a quotient $G_{24}/ \langle \omega
\phi \rangle \cong D_6$, the dihedral group of order 6, we also have
an exact sequence
$$ 0 \rightarrow {\what {\mathbb Z}_3} \rightarrow 
\Ind ^{G_{24}}_{Q_8}{\what {\mathbb Z}_3} \rightarrow 
\Ind ^{G_{24}}_{Q_8}{\what {\mathbb Z}_3}(\theta) \rightarrow 
{\what {\mathbb Z}_3}(\theta) \rightarrow 0,
$$ 
with middle map determined by $1 \otimes 1 \mapsto (c-c^{-1}) \otimes \theta$ 
where $c$ is a generator of $C_3$.

Inducing this to $N$ gives 
\begin{equation} 0 \rightarrow \Ind
^N_{G_{24}}{\what {\mathbb Z}_3} \rightarrow \Ind ^N_{Q_8}{\what
{\mathbb Z}_3} \rightarrow \Ind ^N_{Q_8}{\what {\mathbb Z}_3}(\theta)
\rightarrow \Ind ^N_{G_{24}}{\what {\mathbb Z}_3}(\theta) \rightarrow
0.
\label{eqn:eq2}
\end{equation}

Now splice \eqref{eqn:eq1} and \eqref{eqn:eq2} together at $\Ind
^N_{G_{24}}{\what
{\mathbb Z}_3}(\theta)$ and induce up to ${\mathbb G^1_2}$ to obtain
$$ 0 \rightarrow \Ind ^{\mathbb G^1_2}_{G_{24}}{\what {\mathbb Z}_3}
\rightarrow \Ind ^{\mathbb G^1_2}_{Q_8}{\what {\mathbb Z}_3}
\rightarrow \Ind ^{\mathbb G^1_2}_{Q_8}{\what {\mathbb Z}_3}(\theta)
\rightarrow \Ind ^{\mathbb G^1_2}_{G_{24}}{\what {\mathbb Z}_3}
\rightarrow \Ind ^{\mathbb G^1_2}_N{\what {\mathbb Z}_3} \rightarrow
0.$$
The second and third non-zero terms are projective, so stably 
$\smash{\Ind ^{\mathbb G^1_2}_{G_{24}} {\what {\mathbb Z}_3} \simeq N_3}$. But 
$\smash{\Ind ^{\mathbb G^1_2}_{G_{24}} {\what {\mathbb Z}_3}}$ is cofibrant by
construction and, on restriction to an open torsion free subgroup, 
$N_3$ is a third syzygy hence also cofibrant, so by \fullref{cof} 
there are finitely generated projective modules $P$ and $Q$ such that 
$\smash{\Ind ^{\mathbb G^1_2}_{G_{24}}{\what {\mathbb Z}_3} 
\oplus P \cong N_3 \oplus Q}$.

Let $S$ be a simple $\smash{{\what {\mathbb Z}_3}\llbracket {\mathbb
G^1_2}\rrbracket}$--module
(recall that these correspond to simple $SD_{16}$--modules). Then
$\smash{\Hom _{\mathbb G^1_2}(\Ind ^{\mathbb G^1_2}_{G_{24}}{\what {\mathbb
Z}_3},S)} \cong \smash{\Hom _{G_{24}}({\what {\mathbb Z}_3},S)}$. For this to
be non-zero we need $\smash{\Res ^{\mathbb G^1_2}_{G_{24}}} S \cong
\smash{\mathbb F_3}$, so $S$ must be either ${\mathbb F_3}$ or ${\mathbb
F_3}(\chi )$; in both cases the dimension of the $\Hom$ group is 1.

It follows that the projective cover of $\smash{\Ind ^{\mathbb
G^1_2}_{G_{24}}{\what {\mathbb Z}_3}}$ is $\smash{P_{\mathbb F_3} \oplus
P_{{\mathbb F_3}(\chi )}}$. Now, taking projective covers in $\smash{\Ind
^{\mathbb G^1_2}_{G_{24}}{\what {\mathbb Z}_3}{\oplus}P \cong N_3
\oplus Q}$, we obtain
$$\smash{P_{\mathbb F_3}{\oplus}P_{{\mathbb F_3}(\chi
)}{\oplus}P \cong P_{\mathbb F_3}{\oplus}P_{{\mathbb F_3}(\chi )}
{\oplus}Q},$$
so $P \cong Q$ and thus $\smash{\Ind ^{\mathbb
G^1_2}_{G_{24}}{\what {\mathbb Z}_3}} \cong N_3$, by 
\fullref{ks}.

\begin{remark}
This construction generalizes to $\smash{\mathbb G^1_{p-1}}$ for larger
primes $p$. It is simpler to discuss if we restrict to the Sylow $p$
subgroup. We now have $N=C_p \times \smash{\what {\mathbb Z}_p^{p-2}}$. Since
$\smash{\what {\mathbb Z}_p^{p-2}}$ has cohomological dimension $p-2$, we
could take its projective resolution to the penultimate term and
inflate to $N$. We  then splice on a part induced from a partial
projective resolution of ${\what {\mathbb Z}_3}$ over $C_p$ that is
long enough to make the last term cofibrant. It is not clear whether
this has any significance in homotopy theory.
\end{remark}

\begin{remark}
The Tate--Farrell cohomology of $\mathbb G^1_{p-1}$ is easy to
compute (see Symonds \cite{s1}). It is the low-dimensional cohomology that
is difficult to calculate, but that is precisely what is needed to
identify the projective modules in the complex. If we are satisfied
with a complex with unknown projectives then the construction is
much easier and only depends on the structure of $N$.
\end{remark}

\bibliographystyle{gtart}
\bibliography{link}
\end{document}